\documentclass[twoside]{article}
\usepackage{amsfonts, amsmath}
\usepackage{xcolor}
\usepackage{amssymb}
\usepackage{amsfonts}
\usepackage{textcomp}
\usepackage{microtype}
\usepackage{enumitem}
\usepackage[colorlinks=true]{hyperref}
\hypersetup{linkcolor=blue}
\usepackage[dvips, lmargin=2.5cm, rmargin=1.5cm, tmargin=2.5cm, bmargin=2cm]{geometry}

\newcommand{\verti}[1]{{\left\vert\kern-0.25ex\left\vert\kern-0.25ex\left\vert #1
\right\vert\kern-0.25ex\right\vert\kern-0.25ex\right\vert}}
\newcommand{\vertiii}[1]{{\left\vert\kern-0.25ex\left\vert\kern-0.25ex\left\vert #1
    \right\vert\kern-0.25ex\right\vert\kern-0.25ex\right\vert}}
\title{
\textbf{Three weak solutions for a $(p, q)$-Schr\"{o}dinger-Kirchhoff type
equation}
}
\newenvironment{prv}{ \noindent {\bf Proof } \noindent} {\hfill$\square$ \vskip 6mm}
\newenvironment{th1}{ \noindent {\bf Proof of Theorem \ref{Theorem3}\\} \noindent} {\hfill$\square$ \vskip 6mm}

\numberwithin{equation}{section}
\newtheorem{thm}{Theorem}[section]

\newtheorem{lem}{Lemma}[section]

\newtheorem{cor}{Corolary}[section]
\newtheorem{prop}{Proposition}[section]
\newtheorem{Def}{Definition}[section]
\date{}
\author{Ahmed AHMED$^{*}$$^{1}$, Taghi AHMEDATT$^{2}$ and  Ahmed ABERQI$^{3}$\\
$^{1,2}$University of Nouakchott, Faculty of Science and Technology\\
Mathematics and Computer Sciences Department, \\Research Unit Geometry, Algebra, Analysis  and Applications\\
Nouakchott, Mauritania\\
\href{www.google.com}{\color{blue}ahmedmath2001@gmail.com}$^{*}$$^{1}$  and  \href{http//:www.gmail.com}{\color{blue}taghi-med@hotmail.fr}$^{2}$ \\
$^{3}$National School of Applied Sciences, Sidi Mohamed Ben Abdellah University, Fez, Morocco.\\
\href{http//:www.gmail.com}{\color{blue}ahmed.aberqi@usmba.ac.ma}$^{3}$
}

\everymath{\displaystyle}
\begin{document}
\maketitle
\vspace{0.5cm}
\begin{abstract}
In this manuscript, we investigate a $(p, q)$-Schr\"{o}dinger-Kirchhoff equation involving a continuous
positive potential that meets the del Pino-Felmer type conditions.

Using Recceri's classical variational approach, we prove the existence of three weak solutions.
\end{abstract}
{\bf AMS:}	35Q55,35J60, 35D30,  35J20.
\section{Introduction}
Our objective in this work is to discuss the existence of  three weak solutions for Schr\"{o}dinger double phase Kirchhoff-type elliptic equation involving
potential nonlinearities. The approach is based on a technical strategy that employs the critical points theorem developed by B. Ricceri, which is documented in [\ref{Ricceri29}].

More precisely, we have studied the following equation:
\begin{eqnarray}\label{p}
 &&K_1\Big( \int_{\mathbb{R}^N}\frac{1}{p}|\nabla u|^p dx+ \frac{1}{p}\int_{\mathbb{R}^N} \mathcal{V}(x)|u|^p  dx\Big)\Big(-\Delta_p u
 +\mathcal{V}(x)|u|^{p-2}u\Big)\nonumber\\
 &&+ K_2\Big( \int_{\mathbb{R}^N}\frac{1}{q}|\nabla u|^q dx+ \frac{1}{q}\int_{\mathbb{R}^N} \mathcal{V}(x)|u|^q  dx\Big)\Big(-\Delta_q u
 +\mathcal{V}(x)|u|^{q-2}u\Big)\nonumber\\
 &&\quad= \alpha f_{1}(x,u)+\beta f_{2}(x,u) \quad \textrm{in }\mathbb{R}^N,
\end{eqnarray}
where $1<p<q<N<2q$, $\Delta_su=\mbox{div}(|\nabla u|^{s-2}\nabla u)$ , with $s \in\{p, q\}$, is the s-Laplacian, $\alpha$, $\beta$ are
two real parameters. We suppose that
$\mathcal{V} : \mathbb{R}^{N} \longrightarrow \mathbb{R}$ be a  continuous function that adheres to the conditions elucidated by del Pino and Felmer in their reference [\ref{Pino18}]:
\begin{enumerate}[label=\textup{\textbf{(V\arabic*)}}]
	\item\label{V1} $$\mathcal{V}\in \mathcal{C}(\mathbb{R}^{N}),\ \ \inf_{\mathbb{R}^{N}}\mathcal{V}(x)\geq \mathcal{V}_{0}>0 \textrm{ and }
\mbox{meas}\Big\{x\in\mathbb{R}^{N}:
\mathcal{V}(x)\leq M\Big\}<\infty, \textrm{ for all } M>0,$$
\end{enumerate}
and let $K_1,\; K_2:\mathbb{R}^{+}\longrightarrow \mathbb{R}^{+}$ are a continuous function and satisfies the conditions :
\begin{enumerate}[label=\textup{\textbf{(K\arabic*)}}]
	\item\label{K1}
\begin{equation*}
   k_{1}t^{\sigma_{1}-1} \leq K_{1}(t) \quad \mbox{for all} \quad t>0,
\end{equation*}
where $k_{1}$ and $\sigma_{1}$ are positive constants such that $\sigma_{1}\in \left(1,\frac{q^{*}}{q}\right)$.
\end{enumerate}
\begin{enumerate}[label=\textup{\textbf{(K\arabic*)}}]
\setcounter{enumi}{1}
	\item\label{K2}
\begin{equation*}
   k_{2}t^{\sigma_{2}-1} \leq K_{2}(t)  \quad \mbox{for all} \quad t>0,
\end{equation*}
where $k_{2}$ and $\sigma_{2}$ are positive constants such that  $\sigma_{2}\in \left(1,\frac{q^{*}}{q}\right)$.
\end{enumerate}
\begin{Def}\label{aaa} Let $1 < z < \infty$ with $z < N$. For any Carath\'{e}odory function $f_{1}:\mathbb{R}^{N}\times \mathbb{R}\longrightarrow \mathbb{R}$ such that, we denote the class $\mathcal{A}$ as
\begin{equation*}
   \sup_{(x,t)\in\mathbb{R}^{N}\times \mathbb{R}}\frac{|f_{1}(x, t)|}{|t|^{s-1}}<\infty,
\end{equation*}
for all $s\in [p, q^{*}]$.
\end{Def}
If $f_{1}\in \mathcal{A}$, we assume that:
\begin{enumerate}[label=\textup{\textbf{(F\arabic*)}}]
	\item\label{F1}
\begin{equation*}
   \sup_{u\in W_{\mathcal{V}}^{p,q}(\mathbb{R}^{N})}\int_{\mathbb{R}^{N}}F_{1}(x, u)dx>0,
\end{equation*}
\end{enumerate}
\begin{enumerate}[label=\textup{\textbf{(F\arabic*)}}]
	\setcounter{enumi}{1}
	\item\label{F2}
\begin{equation*}
\limsup_{t\rightarrow0}\frac{\displaystyle\sup_{x\in \mathbb{R}^{N}}F_{1}(x, t)}{|t|^{q\max(\sigma_{1},\sigma_{2})}}=0,
\end{equation*}
\end{enumerate}
\begin{enumerate}[label=\textup{\textbf{(F\arabic*)}}]
	\setcounter{enumi}{2}
	\item\label{F3}
\begin{equation*}
\limsup_{t\rightarrow\infty}\frac{\displaystyle\sup_{x\in \mathbb{R}^{N}}F_{1}(x, t)}{|t|^{p\min(\sigma_{1},\sigma_{2})}}=0.
\end{equation*}
\end{enumerate}

The operator defined in \eqref{p} is used in many branches of mathematics, including calculus of variations and PEDs. It has also been applied in a wide range of physical and engineering contexts, including fluid
filtration in porous media, optimal control, image processing, elastoplasticity,
financial mathematics, and elsewhere, for more details, see [\ref{ddd1}, \ref{ddd2}, \ref{ddd4}] and the references therein.
The Kirchhoff-type problem was primarily introduced in [\ref{Ahmed22}] to generalise the classical D'Alembert wave equation for
free vibrations of elastic strings. Some interesting research by variational methods can be found in [\ref{Azroul}, \ref{ddd7},  \ref{ddd9}, \ref{ddd6}] for
Kirchhoff-type problems. More precisely, Kirchhoff introduced a famous equation defined as
\begin{equation}\label{kirch}
r\frac{\partial^2 u}{\partial t^2}-\left(\frac{P_0}{h}+\frac{E}{2L}\int_0^L\left|\frac{\partial u}{\partial x}\right|^2dx\right)
\frac{\partial^2 u}{\partial x^2}=0,
\end{equation}
that it is related to the problem \eqref{p}. In \eqref{kirch}, $L$ is the length of the string, $h$ is the area of the cross-section, $P_0$ is the initial tension, $r$ is the mass density, and $E$ is
the Young modulus of the material. See the paper [\ref{Ahmed4}, \ref{Ahmed15}, \ref{references1}]
for more details.

In recent research, specifically when considering the functions $K_1(t)=1+at$ and $K_2(t)=1+bt$, alongside the simplifications $f_{1}(x,u)=f(u)$
and $f_{2}(x,u)\equiv0$, where $a$ and $b$ are both strongly positive, W. Zhang, J. Zuo, and P. Zhao, as documented in [\ref{Zhang}], made significant
advancements. They successfully demonstrated the existence of multiple positive solutions to the following problem:
\begin{equation}\label{pppp}
 -\Big(1+a \int_{\mathbb{R}^N}|\nabla u|^p dx\Big)\Delta_p u-\Big(1+b \int_{\mathbb{R}^N}|\nabla u|^q dx\Big)\Delta_q u
 +\mathcal{V}_{\varepsilon}(x)\Big(|u|^{p-2}u+|u|^{q-2}u\Big)\quad=f(u) \quad \textrm{in }\mathbb{R}^N,
\end{equation}
In a similar vein, V. Ambrosio and T. Isernia, as discussed in [\ref{Ambrosio}], also conducted a comprehensive study to establish
the existence of multiple solutions for problem \eqref{pppp}. These groundbreaking contributions have significantly enriched the
understanding of this complex problem.

When $K_1 \equiv K_2 \equiv 1$, the expression in \eqref{p} transforms into a $(p, q)$-Schr\"{o}dinger equation structured as
\begin{equation}\label{p2}
    -\Delta_p u-\Delta_q u+\mathcal{V}(x)\Big(|u|^{p-2}u+|u|^{q-2}u\Big)=\alpha f_{1}(x,u)+\beta f_{2}(x,u) \quad \textrm{in }\mathbb{R}^N.
\end{equation}
It is worth noting that recent years have seen several researchers establish a variety of results related to \eqref{p2}, including existence, multiplicity,
and regularity. In [\ref{doubl28}], some regularity outcomes were established. Furthermore, [\ref{doubl29}] introduced an existence result by combining the
concentration-compactness principle and the mountain pass theorem. Another work by [\ref{doubl22}] demonstrated the existence of a positive ground-state
solution for a $(p, q)$-Schr\"{o}dinger equation with critical growth. Multiple solutions and the concentration of nontrivial solutions for \eqref{p2}
were addressed in [\ref{doubl3}]. The authors in [\ref{doubl39}] employed refined variational methods grounded in critical point theory, along with
 Morse theory and truncation techniques, to establish multiplicity outcomes for a $(p, q)$-Laplacian problem within bounded domains. For additional
 insights, interested readers are encouraged to explore [\ref{doubl7}, \ref{RazaniA},  \ref{aberqidouble}, \ref{references2},  \ref{ouaz}, \ref{oauziz}, \ref{RazaniC}, \ref{RazaniB},  \ref{RazaniD}, \ref{references3}],
 along with the references provided therein.

When $K_1 \equiv K_2 \equiv 1$ and $p=q$, \eqref{p2} boils down a $p$-Schr\"{o}dinger equation of the type
\begin{equation}\label{p3}
    -\Delta_p u+\mathcal{V}(x)\Big(|u|^{p-2}u\Big)=\alpha f_{1}(x,u)+\beta f_{2}(x,u) \quad \textrm{in }\mathbb{R}^N.
\end{equation}
It is important to stress that the world of literature abounds with a variety of captivating and diverse findings.
A comprehensive array of these results is available in sources such as [\ref{doubl23}]. Furthermore, it's essential to
note that when $p=q=2$. Equation \eqref{p3} transforms into the widely studied Schr\"{o}dinger equation, a subject that has garnered extensive
attention over the past three decades. For a deeper exploration, readers are encouraged to examine works like
[\ref{aberqidouble1}, \ref{doubl25}, \ref{doubl43}, \ref{references4}, \ref{references5}].

It is important to note that Numerous studies have employed variational techniques and, by extension, Ricceri's principle to tackle a wide
array of mathematical problems. These methodologies have been applied to Schr\"{o}dinger double phase Kirchhoff-type elliptic equations with
potential nonlinearities, among other equations, to ascertain the existence of weak solutions. The literature demonstrates that variational
Ricceri's principle, along with its associated techniques, offers a robust framework for investigating critical points and establishing the
existence of multiple solutions in a variety of mathematical settings. In this context and in our endeavor of establishing the presence of
weak solutions for our central problem, we draw upon a theorem originally formulated by B. Ricceri in [\ref{Ricceri29}]. This theorem will
serve as a key tool in our approach. For a more in-depth exploration of this theorem, we recommend referring to the following sources:
[\ref{ricceri2000}, \ref{ricceriBasel}, \ref{Ricceri2009}, \ref{Yuan}].

Motivated by the above research, we have successfully demonstrated the existence of three solutions to problem \eqref{p}. What sets this study apart is our comprehensive approach to the Kirchhoff functions, $K_{1}$ and $K_{2}$, which extends beyond the scope of prior research, as exemplified in works like [\ref{Ambrosio}] and [\ref{Zhang}]. Furthermore, our approach yields entirely distinct results. This complexity is further amplified by the inherent lack of compactness in the embedding $W_{\mathcal{V}}^{p,q}(\mathbb{R}^{N})\subset L^{s}(\mathbb{R}^{N})$ for $s\in[p, q^{*}]$ (as elucidated in section 2). These factors collectively present formidable challenges when applying Theorem \ref{ricceri}, which stands as our primary analytical tool. In response to these challenges, we've embraced a fusion of classical and contemporary techniques. This amalgamation becomes especially apparent as we navigate the intricacies, as illustrated in Lemma \ref{ii}, \ref{lem7}, and \ref{lemceorcive}.

Our major result is the following theorems.
\begin{thm}\label{Theorem3}
Assuming that \ref{F1}-\ref{F3}, \ref{V1}, \ref{K1}, and \ref{K2} are satisfied, thus, for any compact interval $\mathbf{\Lambda}\subset(\Theta^{*}, \infty)$, there exists a positive number $\eta>0$ with the following property: for any $\alpha\in \mathbf{\Lambda}$ and any $f_{2}\in \mathcal{A}$, there exists a positive constant $\beta^{*}>0$ such that, for each $\beta\in[0, \beta^{*}]$, equation \eqref{p} admits at least three weak solutions with norms smaller than $\eta$.
\end{thm}

The paper is organized as follows: In Section 2, we collect the main definitions and properties of Sobolev spaces.
In section 3, we give the proof of the main results, and it also includes remarks and illustrative examples.
\section{Some  back  results}
Let  $z$  be  a  real  such  that   $1< z<N<\infty$. We denote by $z'$ the dual exponent of $z$, i.e., $z'=\frac{z}{z-1}$, while $z^{*}$ is the Sobolev dual in dimension $N$, namely
\begin{equation*}
    z^{*}:=\frac{Nz}{N-z}.
\end{equation*}
The term $\mathcal{D}^{1,z}(\mathbb{R}^{N})$ designates the set of functions obtained by taking the closure of $\mathcal{C}^{\infty}_{c} (\mathbb{R}^{N})$ functions, while considering the associated norm
\begin{equation*}
    |\nabla u|^{z}_{z}=\int_{\mathbb{R}^{N}}|\nabla u|^{z}dx.
\end{equation*}
By $W^{1,z}(\mathbb{R}^{N})$, We refer to the Sobolev space, which is equipped with the following norm
\begin{equation*}
    \|u\|_{W^{1,z}(\mathbb{R}^{N})}=(|u|^{z}_{z}+|\nabla u|^{z}_{z})^{\frac{1}{z}}.
\end{equation*}
\begin{thm}\label{Adams} (See [\ref{20222022}]). Assume that $1<z<N$.  $W^{1,z}(\mathbb{R}^{N})$ is continuously embedded in $L^{t}(\mathbb{R}^{N})$
for any $t\in[z, z^{*}]$ and compactly in $L_{loc}^{t}(\mathbb{R}^{N})$ for any $t\in[1, z^{*})$.
\end{thm}

Let $1<p<q<N<2q$ and set
\begin{equation*}
    W^{p,q}(\mathbb{R}^{N})=W^{1,p}(\mathbb{R}^{N})\cap W^{1,q}(\mathbb{R}^{N}),
\end{equation*}
endowed with the norm
\begin{equation*}
    \|u\|_{W^{p,q}(\mathbb{R}^{N})}=\|u\|_{W^{1,p}(\mathbb{R}^{N})}+\|u\|_{W^{1,q}(\mathbb{R}^{N})}.
\end{equation*}

 We introduce the space
\begin{equation*}
    W_{\mathcal{V}}^{p,q}(\mathbb{R}^{N})=\left\{u\in W^{p,q}(\mathbb{R}^{N}): \int_{\mathbb{R}^{N}}\mathcal{V}(x)(|u|^{p}+|u|^{q})dx<\infty\right\}.
\end{equation*}
Characterized by the norm
\begin{equation*}
    \|u\|_{W_{\mathcal{V}}^{p,q}(\mathbb{R}^{N})}=\vertiii{u}_{\mathcal{V},p,q}=\|u\|_{\mathcal{V},p}+\|u\|_{\mathcal{V},q},
\end{equation*}
where
\begin{equation*}
    \|u\|_{\mathcal{V},r}=\left(|\nabla u|^{r}_{r}+\int_{\mathbb{R}^{N}}\mathcal{V}(x)|u|^{r}dx\right)^{\frac{1}{r}} \quad \mbox{for each}\quad  r\in \{p, q\}.
\end{equation*}
\begin{thm}\label{separable} (See [\ref{20222022}]). $W_{\mathcal{V}}^{p,q}(\mathbb{R}^{N})$ is a separable reflexive Banach space.
\end{thm}

Recall that the following embedding results.
\begin{prop}\label{embeddingr}(See [\ref{Zhang}]).
    The space $W_{\mathcal{V}}^{p,q}(\mathbb{R}^{N})$ is continuously embedded into $W^{1,p}(\mathbb{R}^{N})\cap W^{1,q}(\mathbb{R}^{N})$. As a result $W_{\mathcal{V}}^{p,q}(\mathbb{R}^{N})$ is continuously embedded in $L^{s}(\mathbb{R}^{N})$ for any $s\in[p, q^{*}]$
and compactly embedded in $L^{s}_{loc}(\mathbb{R}^{N})$ for any $s\in[1, q^{*})$. Moreover, if $\mathcal{V}_{\infty} =\infty$, the
embedding $W_{\mathcal{V}}^{p,q}(\mathbb{R}^{N})\subset L^{s}(\mathbb{R}^{N})$ is compact for any $p\leq s < q^{*}$.
\end{prop}
\begin{prop}(Algebraic inequality)\label{algebraic}(See [\ref{2626}]).
Let $\alpha \in (0,1)$, for any real numbers $a>0$,$b>0$, $c>0$ and  $d>0$, we have
$$
\left( a+b\right)^{\alpha} \left(c+d\right)^{1-\alpha} \geq a^{\alpha} c^{1-\alpha} + b^{\alpha} d^{1-\alpha}.
$$
\end{prop}
\begin{Def}\label{Definition1}
Consider a real Banach space $X$. We define $\mathcal{W}_{X}$ as the set of functionals $A : X \longrightarrow \mathbb{R}$ with the following characteristic: if $(u_n)_{n}$ is a sequence in $X$ that converges weakly to $u \in X$, and $\liminf_{n\rightarrow\infty}A(u_{n})\leq A(u)$, then there exists a subsequence $(u_{n})_n$ that converges strongly to $u$.
\end{Def}

The following theorem plays an important role in this paper.
\begin{thm}\label{ricceri} (See [\ref{Ricceri29}]).
Consider $X$ as a separable and reflexive real Banach space, with the norm denoted by $\|\cdot\|$. Let $\Psi_{1} : X \longrightarrow \mathbb{R}$ be a coercive, sequentially weakly lower semicontinuous $\mathcal{C}^{1}$ functional that falls under the class $\mathcal{W}_{X}$. Additionally, $\Psi_{1}$ is bounded on each bounded subset of $X$ and has a derivative with a continuous inverse on $X^{*}$. Let $\Psi_{2} : X \longrightarrow \mathbb{R}$ be a $\mathcal{C}^{1}$ functional with a compact derivative.

Suppose that $\Psi_{1}$ possesses a strict local minimum at $x_{0}$, satisfying $\Psi_{1}(x_{0}) = \Psi_{2}(x_{0}) = 0$. Finally, we assume that:
	\begin{equation*}
	\max\left\{\limsup_{\|x\|\rightarrow+\infty}\frac{\Psi_{2}(x)}{\Psi_{1}(x)},\limsup_{x\rightarrow x_{0}}\frac{\Psi_{2}(x)}{\Psi_{1}(x)}\right\}\leq0,
	\end{equation*}
	and that
	\begin{equation*}
	\sup_{x\in X} \min\{\Psi_{1}(x),\Psi_{2}(x)\}>0.
	\end{equation*}
	Let
	\begin{equation*}
	\Theta^{*}:=\inf\left\{\frac{\Psi_{1}(x)}{\Psi_{2}(x)}: x\in X, \min(\Psi_{1}(x),\Psi_{2}(x))>0\right\}.
	\end{equation*}
	Hence, for every compact interval $\mathbf{\Lambda} \subset (\Theta^{*}, +\infty)$, there exists a positive number $\eta$ such that the following holds: for any $\alpha$ within $\mathbf{\Lambda}$ and any $\mathcal{C}^{1}$ functional $\Psi_{3} : X \longrightarrow \mathbb{R}$ with a compact derivative, there is a $\beta^{*}$ greater than zero such that for each $\beta$ in the interval $[0, \beta^{*}]$, the equation
\begin{equation}
\Psi_{1}'(x) = \alpha \Psi_{2}'(x) + \beta \Psi_{3}'(x),
\end{equation}
has at least three solutions with norms smaller than $\eta$.
\end{thm}
\section{Proof of the main results}
We introduce the functionals $\Psi_{1},\Psi_{2}:W_{\mathcal{V}}^{p,q}(\mathbb{R}^{N})\longrightarrow\mathbb{R}$ by
\begin{equation}\label{5aeq1}
\Psi_{1}(u)=\widehat{K}_{1}\left(\varphi_{p,\mathcal{V}}(u)\right)+\widehat{K}_{2}\left(\varphi_{q,\mathcal{V}}(u)\right),
\end{equation}
and
\begin{equation}\label{5aeq2}
\Psi_{2}(u)=\int_{\mathbb{R}^{N}}F_{1}(x,u)dx,
\end{equation}
where
\begin{equation}\label{ppeq3.88}
    \varphi_{p,\mathcal{V}}(u)=\int_{\mathbb{R}^{N}}\frac{1}{p}\Big(|\nabla u|^{p}+\mathcal{V}(x)|u|^{p}\Big)dx,\quad \varphi_{q,\mathcal{V}}(u)=\int_{\mathbb{R}^{N}}\frac{1}{q}\Big(|\nabla u|^{q}+\mathcal{V}(x)|u|^{q}\Big)dx,
\end{equation}
\begin{equation}\label{ppeq3.89}
\widehat{K}_{1}(r)=\int_{0}^{r}K_{1}(s)ds,\quad  \widehat{K}_{2}(r)=\int_{0}^{r}K_{2}(s)ds \quad \forall r\geq0,
\end{equation}
and
\begin{equation}\label{ppeq3.8}
F_{1}(x,t)=\displaystyle{\int^{t}_{0}f_{1}(x,s)ds}\quad  \mbox{for all } (x,t)\in \mathbb{R}^{N}\times\mathbb{R}.
\end{equation}
Under the hypotheses \ref{F1}-\ref{F3}, we set
\begin{equation*}
    \Theta^{*}=\inf\left\{\frac{\Psi_{1}(u)}{\Psi_{2}u)}:u\in W_{\mathcal{V}}^{p,q}(\mathbb{R}^{N}), \int_{\mathbb{R}^{N}}F_{1}(x, u)dx>0\right\}.
\end{equation*}

Let us start by giving the definition of weak solution for the problem \eqref{p}.
\begin{Def}\label{def3.1} We say that $u\in W_{\mathcal{V}}^{p,q}(\mathbb{R}^{N})$ is a weak solution of problem \eqref{p}, if
	\begin{eqnarray*}
		&& K_{1}\left(\int_{\mathbb{R}^{N}}\frac{1}{p}\Big(|\nabla u|^{p}+\mathcal{V}(x)|u|^{p}\Big)dx\right)\int_{\mathbb{R}^{N}}\Big(|\nabla u|^{p-2}\nabla u\nabla v+\mathcal{V}(x)|u|^{p-2}uv\Big)dx \\
		&&+ K_{2}\left(\int_{\mathbb{R}^{N}}\frac{1}{q}\Big(|\nabla u|^{q}+\mathcal{V}(x)|u|^{q}\Big)dx\right)\int_{\mathbb{R}^{N}}\Big(|\nabla u|^{q-2}\nabla u\nabla v+\mathcal{V}(x)|u|^{q-2}uv\Big)dx\\
		&&\quad=\alpha \displaystyle\int_{\Omega}f_{1}(x, u)v\>dx+\beta \displaystyle\int_{\Omega}f_{2}(x, u)v\>dx,
	\end{eqnarray*}
	for all $v\in W_{\mathcal{V}}^{p,q}(\mathbb{R}^{N})$.
\end{Def}
The proof of Theorem \ref{Theorem3} is methodically divided into three sequential steps.

\begin{th1}
\begin{center}
\bf Step 1 : Technical lemmas
\end{center}
This subsection is dedicated to laying the groundwork by introducing essential technical lemmas that will find repeated application in this paper.

\begin{lem} \label{alem3.1}
	Let $f_{1}\in \mathcal{A}$, then the fonctional $\Psi_{2}\in \mathcal{C}^{1}(W_{\mathcal{V}}^{p,q}(\mathbb{R}^{N}),\mathbb{R})$ with the derivative is given by
	\begin{equation}\label{equ }
\langle \Psi_{2}'(u), v\rangle=\displaystyle\int_{\Omega}f_{1}(x, u)v\>dx,
	\end{equation}
	for all $u,v\in W_{\mathcal{V}}^{p,q}(\mathbb{R}^{N})$. In addition  $\Psi_{2}':W_{\mathcal{V}}^{p,q}(\mathbb{R}^{N})\longrightarrow\Big(W_{\mathcal{V}}^{p,q}(\mathbb{R}^{N})\Big)^{*}$ is compact.
\end{lem}
\begin{prv}
First, thanks  to  the  definition of $\mathcal{A}$ and  the  embedding theorem, $\Psi_{2}$ is  well  difined on $W_{\mathcal{V}}^{p,q}(\mathbb{R}^{N})$.
New  we show  that  $\Psi_{2}$ is G\^{a}teaux-differentiable on $W_{\mathcal{V}}^{p,q}(\mathbb{R}^{N})$ with the
derivative is given by (\ref{equ }). Actually, let $(u_n)_{n}\subset W_{\mathcal{V}}^{p,q}(\mathbb{R}^{N})$ be a sequence converging strongly to $u\in W_{\mathcal{V}}^{p,q}(\mathbb{R}^{N})$. Because $W_{\mathcal{V}}^{p,q}(\mathbb{R}^{N})$ is compactly embedded in $L^{s}(\mathbb{R}^N),$  thus $(u_n)_{n}$ converges strongly
to $u$ in  $L^s(\mathbb{R}^N)$. Consequently, there exists a subsequence of $(u_n)_{n}$, denoted as $(u_n)_{n}$, along with a function $\overline{u}\in L^{s}(\mathbb{R}^N)$. This subsequence converges to $\overline{u}$ almost everywhere in $\mathbb{R}^N$, satisfying $|u_n|\leq |\overline{u}|$ for every natural number $n$ and nearly everywhere in $\mathbb{R}^N$.  Given $f_{1}\in \mathcal{A}$, for every measurable function $u : \mathbb{R}^N \longrightarrow \mathbb{R}$, the operator defined as $u \longmapsto f_{1}(\cdot, u(\cdot))$ maps $L^{s}(\mathbb{R}^N)$ into $L^{s'}(\mathbb{R}^N)$. Fix $v\in W_{\mathcal{V}}^{p,q}(\mathbb{R}^{N})$ with $\|v\|_{W_{\mathcal{V}}^{p,q}(\mathbb{R}^{N})} \leq 1$, using H\"{o}lder's inequality and the continuous embedding of $W_{\mathcal{V}}^{p,q}(\mathbb{R}^{N})$ into  $L^{s}(\mathbb{R}^N)$ , we have
 \begin{eqnarray*}
 	\Big|\langle \Psi_{2}'(u_n)-\Psi_{2}'(u), v\rangle\Big| &=& \Big|\int_{\mathbb{R}^{N}}\Big(f_1(x,u_n(x))-f_1(x,u(x) )\Big) v(x)dx\Big| \\
 	&\leq & \|f_1(x,u_n(x))-f_1(x,u(x) )\|_{L^{s'}(\mathbb{R}^{N})}\|v\|_{L^{s}(\mathbb{R}^{N})}\\
 	&\leq & C_1\|f_1(x,u_n(x))-f_1(x,u(x) )\|_{L^{s'}(\mathbb{R}^{N})}\|v\|_{W_{\mathcal{V}}^{p,q}(\mathbb{R}^{N})},
 \end{eqnarray*}
a certain positive constant $C_1$, we can pass to the supremum for $\|v\|_{W_{\mathcal{V}}^{p,q}(\mathbb{R}^{N})} \leq 1$, we get
 $$\|\langle \Psi_{2}'(u_n)-\Psi_{2}'(u), v\rangle\|_{\big(W_{\mathcal{V}}^{p,q}(\mathbb{R}^{N})\big)^{*}}\leq\|f_1(x,u_n(x))-f_1(x,u(x) )\|_{L^{s'}(\mathbb{R}^{N})}.$$
 Since  $f_{1}\in \mathcal{A}$, we deduce
 $$ f_1(x,u_n(x))-f_1(x,u(x) ) \longrightarrow  0 \text{ as } n\longrightarrow \infty ,$$ and
 $$\Big|f_1(x,u_n(x))-f_1(x,u(x) )\Big|\leq C_2\Big(|\overline{u}(x)|^{s-1}+|u(x)|^{s-1}\Big).$$
Considering almost every $x\in \mathbb{R}^{N}$, it is worth noting that the majorant function in the previous relation lies in $L^{s'}(\mathbb{R}^N)$. Consequently, utilizing the dominated convergence theorem, we obtain that
 $$\|f_1(x,u_n(x))-f_1(x,u(x) )\|_{s'}\longrightarrow 0  \text{ as } n\longrightarrow \infty.$$
 The continuity of $\Psi_{2}'$ is demonstrated. To establish the compactness of $\Psi_{2}'$, consider a sequence $(u_n)_{n}\subset W_{\mathcal{V}}^{p,q}(\mathbb{R}^{N})$ that is bounded. Then there exists a subsequence of $(u_n)_{n}$ (still denoted
 by $(u_n)_{n}$ ) weakly convergent in $W_{\mathcal{V}}^{p,q}(\mathbb{R}^{N})$ and strongly convergent in $L^{s}(\mathbb{R}^N)$. Using a similar argument as outlined above, we infer that the sequence $(\Psi_{2}'(u_n))_{n}$ exhibits strong convergence, and the compactness of the operator $\Psi_{2}'$ is established.
\end{prv}
\begin{lem}\label{alem3.2} Let \ref{V1}, \ref{K1} and \ref{K2} hold true. Then the functional $\Psi_{1} \in  \mathcal{C}^{1}(W_{\mathcal{V}}^{p,q}(\mathbb{R}^{N}),\mathbb{R})$ and
\begin{eqnarray}\label{a5a5}
  \langle \Psi_{1}'(u), v\rangle &=& K_{1}\left(\int_{\mathbb{R}^{N}}\frac{1}{p}\Big(|\nabla u|^{p}+\mathcal{V}(x)|u|^{p}\Big)dx\right)\int_{\mathbb{R}^{N}}\Big(|\nabla u|^{p-2}\nabla u\nabla v+\mathcal{V}(x)|u|^{p-2}uv\Big)dx \nonumber\\
  &+&K_{2}\left(\int_{\mathbb{R}^{N}}\frac{1}{q}\Big(|\nabla u|^{q}+\mathcal{V}(x)|u|^{q}\Big)dx\right)\int_{\mathbb{R}^{N}}\Big(|\nabla u|^{q-2}\nabla u\nabla v+\mathcal{V}(x)|u|^{q-2}uv\Big)dx,
\end{eqnarray}
 for each $u, v\in W_{\mathcal{V}}^{p,q}(\mathbb{R}^{N})$.  Morevoer,  for all $u\in W_{\mathcal{V}}^{p,q}(\mathbb{R}^{N})$, $\Psi_{1}'\in \Big(W_{\mathcal{V}}^{p,q}(\mathbb{R}^{N})\Big)^{*}$.
\end{lem}
\begin{prv}
Initially, it is straightforward to observe that
\begin{eqnarray*}
  \langle \Psi_{1}'(u), v\rangle &=& K_{1}\left(\int_{\mathbb{R}^{N}}\frac{1}{p}\Big(|\nabla u|^{p}+\mathcal{V}(x)|u|^{p}\Big)dx\right)\int_{\mathbb{R}^{N}}\Big(|\nabla u|^{p-2}\nabla u\nabla v+\mathcal{V}(x)|u|^{p-2}uv\Big)dx \\
  &+&K_{2}\left(\int_{\mathbb{R}^{N}}\frac{1}{q}\Big(|\nabla u|^{q}+\mathcal{V}(x)|u|^{q}\Big)dx\right)\int_{\mathbb{R}^{N}}\Big(|\nabla u|^{q-2}\nabla u\nabla v+\mathcal{V}(x)|u|^{q-2}uv\Big)dx,
\end{eqnarray*}
for all $u, v\in W_{\mathcal{V}}^{p,q}(\mathbb{R}^{N})$. It follows from \eqref{a5a5} that for each $u\in W_{\mathcal{V}}^{p,q}(\mathbb{R}^{N})$, $\Psi_{1}'\in \Big(W_{\mathcal{V}}^{p,q}(\mathbb{R}^{N})\Big)^{*}$.

Next, we prove that $\Psi_{1}':W_{\mathcal{V}}^{p,q}(\mathbb{R}^{N})\longrightarrow \left(W_{\mathcal{V}}^{p,q}(\mathbb{R}^{N})\right)^{*}$ is continuous. Let $(u_{n})_{n}\in W_{\mathcal{V}}^{p,q}(\mathbb{R}^{N})$ with $u_{n}\longrightarrow u$ strongly in $W_{\mathcal{V}}^{p,q}(\mathbb{R}^{N})$.
We can assume, without sacrificing generality, that $u_{n}$ converges to $u$ almost everywhere in $\mathbb{R}^{N}$. Consequently, the sequences
\begin{equation}\label{sequencesp}
    (|\nabla u_{n}|^{p-2}\nabla u_{n})_{n} \mbox{ is bounded in } L^{p^{'}}(\mathbb{R}^{N}\times\mathbb{R}^{N}),
\end{equation}
and
\begin{equation}\label{sequencesq}
   (|\nabla u_{n}|^{q-2}\nabla u_{n})_{n} \mbox{ is bounded in } L^{q^{'}}(\mathbb{R}^{N}\times\mathbb{R}^{N}),
\end{equation}
as well as a.e. in $\mathbb{R}^{N}\times \mathbb{R}^{N}$,
\begin{equation*}
    (|\nabla u_{n}|^{p-2}\nabla u_{n})\longrightarrow (|\nabla u|^{p-2}\nabla u),
\end{equation*}
and
\begin{equation*}
    (|\nabla u_{n}|^{q-2}\nabla u_{n})\longrightarrow(|\nabla u|^{q-2}\nabla u).
\end{equation*}
Thus, the Brezis-Lieb lemma (see [\ref{Brezis}]) implies
\begin{equation*}
   \int_{\mathbb{R}^{N}} \lim_{n\longrightarrow+\infty}\Big|(|\nabla u_{n}|^{p-2}\nabla u_{n}-|\nabla u|^{p-2}\nabla u)\Big|^{p'}dx =\lim_{n\rightarrow+\infty} \int_{\mathbb{R}^{N}} (|\nabla u_{n}|^{p}-|\nabla u|^{p})dx,
\end{equation*}
and
\begin{equation*}
    \int_{\mathbb{R}^{N}}\lim_{n\longrightarrow+\infty}\Big|(|\nabla u_{n}|^{q-2}\nabla u_{n}-|\nabla u|^{q-2}\nabla u) \Big|^{q'}dx =\lim_{n\rightarrow+\infty} \int_{\mathbb{R}^{N}} (|\nabla u_{n}|^{q}-|\nabla u|^{q})dx.
\end{equation*}
The strong convergence $u_{n}\longrightarrow u$ in $W_{\mathcal{V}}^{p,q}(\mathbb{R}^{N})$ implies that
\begin{equation}\label{a6a6}
    \lim_{n\rightarrow+\infty} \int_{\mathbb{R}^{N}} (|\nabla u_{n}|^{p}-|\nabla u|^{p})dx=0,
\end{equation}
and
\begin{equation}\label{a6a6a6}
    \lim_{n\rightarrow+\infty} \int_{\mathbb{R}^{N}} (|\nabla u_{n}|^{q}-|\nabla u|^{q})dx=0.
\end{equation}
In a similar fashion,
\begin{equation}\label{a7a7}
    \lim_{n\rightarrow+\infty}\int_{\mathbb{R}^{N}} \mathcal{V}(x)\Big||u_{n}(x)|^{p-2}u_{n}(x)-|u(x)|^{p-2}u(x)\Big|^{p'}dx=0,
\end{equation}
and
\begin{equation}\label{a77}
    \lim_{n\rightarrow+\infty}\int_{\mathbb{R}^{N}} \mathcal{V}(x)\Big||u_{n}(x)|^{q-2}u_{n}(x)-|u(x)|^{q-2}u(x)\Big|^{q'}dx=0.
\end{equation}
Moreover, the continuity of $K_{1}$ and $K_{2}$ implies that
\begin{equation}\label{q777}
    \lim_{n\rightarrow+\infty}K_{1}\left(\int_{\mathbb{R}^{N}}\frac{1}{p}(|\nabla u_{n}|^{p}+\mathcal{V}(x)|u_{n}|^{p})dx\right)
    \longrightarrow K_{1}\left(\int_{\mathbb{R}^{N}}\frac{1}{p}(|\nabla u|^{p}+\mathcal{V}(x)|u|^{p})dx\right),
\end{equation}
and
\begin{equation}\label{7777}
     \lim_{n\rightarrow+\infty}K_{2}\left(\int_{\mathbb{R}^{N}}\frac{1}{q}(|\nabla u_{n}|^{q}+\mathcal{V}(x)|u_{n}|^{q})dx\right)
     \longrightarrow K_{2}\left(\int_{\mathbb{R}^{N}}\frac{1}{q}(|\nabla u|^{q}+\mathcal{V}(x)|u|^{q})dx\right).
\end{equation}
Utilizing H\"{o}lder's inequality along with the combination of \eqref{a6a6}-\eqref{7777}, we get
\begin{equation*}
    \|\Psi_{1}(u_{n})-\Psi_{1}(u)\|_{\left(W_{\mathcal{V}}^{p,q}(\mathbb{R}^{N})\right)^{*}}=\sup_{u\in W_{\mathcal{V}}^{p,q}(\mathbb{R}^{N}), \vertiii{v}_{\mathcal{V},p,q}\leq1 }| \langle \Psi_{1}'(u_{n})-\Psi_{1}'(u), v\rangle|\longrightarrow0.
\end{equation*}

This completes the proof that $\Psi_{1}':W_{\mathcal{V}}^{p,q}(\mathbb{R}^{N})\longrightarrow \left(W_{\mathcal{V}}^{p,q}(\mathbb{R}^{N})\right)^{*}$ is continuous, and therefore $\Psi_{1} \in  \mathcal{C}^{1}(W_{\mathcal{V}}^{p,q}(\mathbb{R}^{N}),\mathbb{R})$.
\end{prv}
\begin{lem}\label{alem3.2} The functional $\Psi_{1}$ is sequentially weakly lower semicontinuous.
\end{lem}
\begin{prv}
Consider a sequence $(u_{n})_{n}$ such that $u_{n}\rightharpoonup u$ in $W_{\mathcal{V}}^{p,q}(\mathbb{R}^{N})$. Due to the convexity of $\varphi_{p,\mathcal{V}}$ and $\varphi_{q,\mathcal{V}}$, for every $n$, it follows that
$$\varphi_{p,\mathcal{V}}(u)\leq \varphi_{p,\mathcal{V}}(u_{n})+\langle \varphi_{p,\mathcal{V}}'(u), u-u_{n}\rangle,$$
and
$$\varphi_{q,\mathcal{V}}(u)\leq \varphi_{q,\mathcal{V}}(u_{n})+\langle \varphi_{q,\mathcal{V}}'(u), u-u_{n}\rangle.$$
Taking the limit in the inequality above as $n$ tends to infinity, we observe that $\varphi_{p,\mathcal{V}}$ and $\varphi_{q,\mathcal{V}}$ are sequentially weakly lower semi-continuous. Then we have
	\begin{equation}\label{sem3.8.3.8}
	\int_{\mathbb{R}^{N}}\frac{1}{p}|\nabla u|^{p}dx\leq \liminf_{n\rightarrow+\infty} \int_{\Omega}\frac{1}{p}|\nabla u_{n}|^{p}dx,
	\end{equation}
\begin{equation}\label{sem3.8.3.9}
	\int_{\mathbb{R}^{N}}\frac{\mathcal{V}(x)}{p}|u|^{p}dx\leq \liminf_{n\rightarrow+\infty} \int_{\Omega}\frac{\mathcal{V}(x)}{p}|u_{n}|^{p}dx,
	\end{equation}
\begin{equation}\label{sem3.8.3.10}
	\int_{\mathbb{R}^{N}}\frac{1}{q}|\nabla u|^{q}dx\leq \liminf_{n\rightarrow+\infty} \int_{\Omega}\frac{1}{q}|\nabla u_{n}|^{q}dx,
	\end{equation}
and
\begin{equation}\label{sem3.8.3.11}
	\int_{\mathbb{R}^{N}}\frac{\mathcal{V}(x)}{q}|u|^{q}dx\leq \liminf_{n\rightarrow+\infty} \int_{\Omega}\frac{\mathcal{V}(x)}{q}|u_{n}|^{q}dx.
	\end{equation}
	By \eqref{sem3.8.3.8}-\eqref{sem3.8.3.11} and since $\widehat{K}_{1}$ and $\widehat{K}_{2}$ are continuous and monotone, we have
	\begin{eqnarray}\label{sem3.88.3.88}
	\liminf_{n\rightarrow+\infty} \Psi_{1}(u_{n})&=&\liminf_{n\rightarrow+\infty} \widehat{K}_{1}\left(\int_{\mathbb{R}^{N}}\frac{1}{p}\Big(|\nabla u_{n}|^{p}+\mathcal{V}(x)|u_{n}|^{p}\Big)dx\right)+\liminf_{n\rightarrow+\infty}\widehat{K}_{2}\left(\int_{\mathbb{R}^{N}}\frac{1}{q}\Big(|\nabla u_{n}|^{q}+\mathcal{V}(x)|u_{n}|^{q}\Big)dx\right)\nonumber\\
&\geq&\widehat{K}_{1}\left(\liminf_{n\rightarrow+\infty}\int_{\mathbb{R}^{N}}\frac{1}{p}\Big(|\nabla u_{n}|^{p}+\mathcal{V}(x)|u_{n}|^{p}\Big)dx\right)+\widehat{K}_{2}\left(\liminf_{n\rightarrow+\infty}\int_{\mathbb{R}^{N}}\frac{1}{q}\Big(|\nabla u_{n}|^{q}+\mathcal{V}(x)|u_{n}|^{q}\Big)dx\right)\nonumber\\
&\geq&\widehat{K}_{1}\left(\liminf_{n\rightarrow+\infty}\int_{\mathbb{R}^{N}}\frac{1}{p}|\nabla u_{n}|^{p}dx +\liminf_{n\rightarrow+\infty}\int_{\mathbb{R}^{N}}\frac{\mathcal{V}(x)}{p}|u_{n}|^{p}dx\right)\nonumber\\
&&+\widehat{K}_{2}\left(\liminf_{n\rightarrow+\infty}\int_{\mathbb{R}^{N}}\frac{1}{q}|\nabla u_{n}|^{q}dx +\liminf_{n\rightarrow+\infty}\int_{\mathbb{R}^{N}}\frac{\mathcal{V}(x)}{q}|u_{n}|^{q}dx\right)\nonumber\\
&\geq&\widehat{K}_{1}\left(\int_{\mathbb{R}^{N}}\frac{1}{p}|\nabla u|^{p}dx +\int_{\mathbb{R}^{N}}\frac{\mathcal{V}(x)}{p}|u|^{p}dx\right)+\widehat{K}_{2}\left(\int_{\mathbb{R}^{N}}\frac{1}{q}|\nabla u|^{q}dx +\int_{\mathbb{R}^{N}}\frac{\mathcal{V}(x)}{q}|u|^{q}dx\right)\nonumber\\
&=& \widehat{K}_{1}\left(\int_{\mathbb{R}^{N}}\frac{1}{p}\Big(|\nabla u|^{p}+\mathcal{V}(x)|u|^{p}\Big)dx\right)+\widehat{K}_{2}\left(\int_{\mathbb{R}^{N}}\frac{1}{q}\Big(|\nabla u|^{q}+\mathcal{V}(x)|u|^{q}\Big)dx\right),\nonumber\\
&=&\widehat{K}_{1}\left(\varphi_{p,\mathcal{V}}(u)\right)+\widehat{K}_{2}\left(\varphi_{q,\mathcal{V}}(u)\right)\nonumber\\
&=& \Psi_{1}(u),
\end{eqnarray}
	namely, $\Psi_{1}$ is sequentially weakly lower semi-continuous.
\end{prv}

In the following lemma, we will demonstrate that the functional $\Psi_{1}$ belongs to the class $\mathcal{W}_{X}$, where
$X =W_{\mathcal{V}}^{p,q}(\mathbb{R}^{N})$ as defined in Definition \ref{Definition1}.
\begin{lem}\label{ii}
The functional $\Psi_{1}$ belongs to the class $\mathcal{W}_{W_{\mathcal{V}}^{p,q}(\mathbb{R}^{N})}$.
\end{lem}
\begin{prv}
Since $\widehat{K}_{1}$ and $\widehat{K}_{2}$ are continuous and strictly increasing, it suffices to show that $\varphi_{p,\mathcal{V}}, \varphi_{q,\mathcal{V}}\in \mathcal{W}_{W_{\mathcal{V}}^{p,q}(\mathbb{R}^{N})}$.
Consider a sequence $(u_n)_{n}$ that converges weakly to $u$ within the space $W_{\mathcal{V}}^{p,q}(\mathbb{R}^{N})$ and let $\liminf_{n\rightarrow+\infty}\varphi_{p,\mathcal{V}}(u_{n})\leq\varphi_{p,\mathcal{V}}(u)$ and $\liminf_{n\rightarrow+\infty}\varphi_{q,\mathcal{V}}(u_{n})\leq\varphi_{q,\mathcal{V}}(u)$. Since the functional $\varphi_{p,\mathcal{V}}$, $ \varphi_{q,\mathcal{V}}$ is sequentially weakly lower semicontinuous, then there exists a subsequence of $(u_n)_{n}$ , bearing the same notation $(u_n)_{n}$ such that $\liminf_{n\rightarrow+\infty}\varphi_{p,\mathcal{V}}(u_{n})\geq\varphi_{p,\mathcal{V}}(u)$ and $\liminf_{n\rightarrow+\infty}\varphi_{q,\mathcal{V}}(u_{n})\geq\varphi_{q,\mathcal{V}}(u)$,
so:
\begin{equation}\label{147}
    \liminf_{n\rightarrow+\infty}\varphi_{p,\mathcal{V}}(u_{n})=\varphi_{p,\mathcal{V}}(u),
\end{equation}
and
\begin{equation}\label{148}
    \liminf_{n\rightarrow+\infty}\varphi_{q,\mathcal{V}}(u_{n})=\varphi_{q,\mathcal{V}}(u).
\end{equation}
According to \eqref{147}-\eqref{148}, we obtain
\begin{equation*}
    \lim_{n\rightarrow\infty}\Psi_{1}(u_{n})=\Psi_{1}(u),
\end{equation*}
that is,
\begin{equation*}
   \lim_{n\rightarrow\infty} \vertiii{u_{n}}_{\mathcal{V},p,q}=   \vertiii{u}_{\mathcal{V},p,q}.
\end{equation*}
Since $W_{\mathcal{V}}^{p,q}(\mathbb{R}^{N})$ is uniformly convex of the space, we can deduce the strong convergence of $u_{n}$ to $u$ in $W_{\mathcal{V}}^{p,q}(\mathbb{R}^{N})$.
\end{prv}
\begin{lem}\label{lem7}
\begin{description}
  \item[(a)] Let $u, v \in W_{\mathcal{V}}^{p,q}(\mathbb{R}^{N})$, this implies that
  \begin{equation}\label{10}
    \langle \varphi'_{p,\mathcal{V}}(u)-\varphi'_{p,\mathcal{V}}(u)'(v), u-v\rangle\geq\Big(\|u\|_{\mathcal{V},p}^{p-1}-\|u\|_{\mathcal{V},p}^{p-1}\Big)\Big(\|u\|_{\mathcal{V},p}-\|u\|_{\mathcal{V},p}\Big).
  \end{equation}
  \begin{equation}\label{10bis}
    \langle \varphi'_{q,\mathcal{V}}(u)-\varphi'_{q,\mathcal{V}}(u)'(v), u-v\rangle\geq\Big(\|u\|_{\mathcal{V},q}^{q-1}-\|u\|_{\mathcal{V},q}^{q-1}\Big)\Big(\|u\|_{\mathcal{V},q}-\|u\|_{\mathcal{V},q}\Big).
  \end{equation}
  \item[(b)] Suppose the sequence $(u_n)_{n}$ weakly converges to $u$ in $W_{\mathcal{V}}^{p,q}(\mathbb{R}^{N})$, and
  \begin{equation}\label{9}
    \limsup_{n\rightarrow\infty}\langle \varphi'_{p,\mathcal{V}}(u_{n}), u_{n}-u\rangle\leq0.
  \end{equation}
   \begin{equation}\label{9bis}
    \limsup_{n\rightarrow\infty}\langle \varphi'_{q,\mathcal{V}}(u_{n}), u_{n}-u\rangle\leq0.
  \end{equation}
As a result, the sequence $(u_n)_{n}$ strongly converges to $u$ in $W_{\mathcal{V}}^{p,q}(\mathbb{R}^{N})$.
\end{description}
\end{lem}
\begin{prv}
  \begin{itemize}
      \item[(a)]  For $u, v \in W_{\mathcal{V}}^{p,q}(\mathbb{R}^{N})$, we have
   \begin{align*}
   &\langle \varphi'_{p,\mathcal{V}}(u)-\varphi'_{p,\mathcal{V}}(u)'(v), u-v\rangle\\
  &=\int_{\mathbb{R}^N} \left( \vert \nabla u \vert^{p-2}\nabla u + \mathcal{V}(x)\vert u\vert^{p-2}u \right) \left(u-v\right)\,dx-\int_{\mathbb{R}^N} \left( \vert \nabla v \vert^{p-2}\nabla v + \mathcal{V}(x)\vert v\vert^{p-2}v \right) \left(u-v\right)\,dx\\
 &= \|u\|_{\mathcal{V},p}^{p} +\|u\|_{\mathcal{V},p}^{p} -\int_{\mathbb{R}^N}  \vert \nabla u \vert^{p-2}\nabla u .v \,dx-\int_{\mathbb{R}^N}  \vert \nabla v \vert^{p-2}\nabla v .u \,dx-\int_{\mathbb{R}^N} \mathcal{V}(x)\vert u\vert^{p-2}u v\,dx-\int_{\mathbb{R}^N} \mathcal{V}(x)\vert v\vert^{p-2}v u \,dx.
  \end{align*}
 The algebraic inequality Proposition ~\ref{algebraic}, with $\alpha=\frac{p-1}{p}$ combining with H\"older inequality, allows us to deduce
 $$\int_{\mathbb{R}^N}  \vert \nabla u \vert^{p-2}\nabla u .v \,dx +\int_{\mathbb{R}^N} \mathcal{V}(x)\vert u\vert^{p-2}u v\,dx \leq \|u\|_{\mathcal{V},p}^{p-1} \|v\|_{\mathcal{V},p},$$
 and
 $$\int_{\mathbb{R}^N}  \vert \nabla v \vert^{p-2}\nabla v .u \,dx +\int_{\mathbb{R}^N} \mathcal{V}(x)\vert v\vert^{p-2}v. u\,dx \leq \|v\|_{\mathcal{V},p}^{p-1} \|u\|_{\mathcal{V},p}.$$
 Then we get~\eqref{10}. Following a similar line of thought, we get ~\eqref{10bis}.
 \item[(b)] Let  $(u_n)_{n}$ converges weakly to $u$ in $W_{\mathcal{V}}^{p,q}(\mathbb{R}^{N})$ such that
  \begin{equation}\label{lims}
    \limsup_{n\rightarrow\infty}\langle \varphi'_{p,\mathcal{V}}(u_{n}), u_{n}-u\rangle\leq0, \,\text{and}\,
  \limsup_{n\rightarrow\infty}\langle \varphi'_{q,\mathcal{V}}(u_{n}), u_{n}-u\rangle\leq 0.
  \end{equation}
  Since $W_{\mathcal{V}}^{p,q}(\mathbb{R}^{N})$  is reflexive Banach space, hence its  isomorphic
to a locally uniformly convex space. Then, all you need to prove the norm convergence, which is attainable by using the assertion (a) above with~\eqref{lims}.
  \end{itemize}
\end{prv}

In this lemma, we demonstrate the invertibility and continuity of the energy $\Psi_{1}'$.
\begin{lem}\label{invertible}
Let \ref{V1}, \ref{K1} and \ref{K2} hold true. Then the operator $\Psi_{1}': W_{\mathcal{V}}^{p,q}(\mathbb{R}^{N})\longrightarrow \Big(W_{\mathcal{V}}^{p,q}(\mathbb{R}^{N})\Big)^{*}$ is, in fact, invertible, and ${\Psi_{1}'}^{-1}$ is notably continuous.
continuous.
\end{lem}
\begin{prv}
To begin, we make the initial assumption that the operator $\Psi_{1}': W_{\mathcal{V}}^{p,q}(\mathbb{R}^{N})\longrightarrow \Big(W_{\mathcal{V}}^{p,q}(\mathbb{R}^{N})\Big)^{*}$ possesses invertibility within the space $W_{\mathcal{V}}^{p,q}(\mathbb{R}^{N})$. As stipulated by the Minty-Browder Theorem (refer to [\ref{Brezis}, \ref{Zeidler}]), our task is to demonstrate that $\Psi_{1}'$ fulfills the criteria of being strictly monotone, hemicontinuous, and coercive in the context of monotone operators. With this in mind, consider $u$ and $v$ as elements of $W_{\mathcal{V}}^{p,q}(\mathbb{R}^{N})$, where $u \neq v$, and let $\mu_{1}$ and $\mu_{2}$ be in the interval $[0, 1]$, with $\mu_{1} + \mu_{2} = 1$. Using \eqref{10} and \eqref{10bis} , we have
$\varphi'_{p,\mathcal{V}}, \varphi'_{q,\mathcal{V}}(u): W_{\mathcal{V}}^{p,q}(\mathbb{R}^{N})\longrightarrow \Big(W_{\mathcal{V}}^{p,q}(\mathbb{R}^{N})\Big)^{*}$ are strictly monotone, so by [\ref{Zeidler}], $f$ is strictly
convex. Moreover, since $K_{1}$ and $K_{2}$ are nondecreasing, the functions $\widehat{K}_{1}$ and $\widehat{K}_{2}$ are convex in $\mathbb{R}^{+}$, thus
we have
\begin{eqnarray}\label{strictlyconvex}
  \Psi_{1}(\mu_{1} u+\mu_{2} v)&=&\widehat{K}_{1}\left(\varphi_{p,\mathcal{V}}(\mu_{1} u+\mu_{2} v)\right)+\widehat{K}_{2}\left(\varphi_{q,\mathcal{V}}(\mu_{1} u+\mu_{2} v)\right)\nonumber\\
  &<& \widehat{K}_{1}\left(\mu_{1}\varphi_{p,\mathcal{V}}( u)+\mu_{2}\varphi_{p,\mathcal{V}}( v)\right)+\widehat{K}_{2}\left(\mu_{1} \varphi_{q,\mathcal{V}}( u)+\mu_{2}\varphi_{q,\mathcal{V}}( v)\right)\nonumber\\
  &\leq& \mu_{1} \widehat{K}_{1}\left(\varphi_{p,\mathcal{V}}( u)\right)+\mu_{2}\widehat{K}_{1}\left(\varphi_{p,\mathcal{V}}( v)\right)+
  \mu_{1} \widehat{K}_{2}\left(\varphi_{q,\mathcal{V}}( u)\right)+\mu_{2}\widehat{K}_{2}\left(\varphi_{q,\mathcal{V}}( v)\right)\nonumber\\
  &\leq& \mu_{1}\left( \widehat{K}_{1}\left(\varphi_{p,\mathcal{V}}( u)\right)+\widehat{K}_{2}\left(\varphi_{q,\mathcal{V}}( u)\right)\right)+
  \mu_{2}\left( \widehat{K}_{1}\left(\varphi_{p,\mathcal{V}}( v)\right)+\widehat{K}_{2}\left(\varphi_{q,\mathcal{V}}( v)\right)\right)\nonumber\\
  &\leq& \mu_{1}\Psi_{1}(u)+\mu_{2}\Psi_{1}(v).
\end{eqnarray}
This shows that $\Psi_{1}$ is strictly convex and is already said, that $\Psi_{1}'$ is strictly monotone.

Let $\|u\|_{\mathcal{V},p}>1$ and $\|u\|_{\mathcal{V},q}>1$, by \ref{K1} and \ref{K2}, we have
\begin{eqnarray}\label{coer}
\langle \Psi_{1}'(u), u\rangle &=& K_{1}\left(\int_{\mathbb{R}^{N}}\frac{1}{p}\Big(|\nabla u|^{p}+\mathcal{V}(x)|u|^{p}\Big)dx\right)\int_{\mathbb{R}^{N}}\Big(|\nabla u|^{p-2}\nabla u\nabla u+\mathcal{V}(x)|u|^{p-2}uu\Big)dx \nonumber\\
  &+&K_{2}\left(\int_{\mathbb{R}^{N}}\frac{1}{q}\Big(|\nabla u|^{q}+\mathcal{V}(x)|u|^{q}\Big)dx\right)\int_{\mathbb{R}^{N}}\Big(|\nabla u|^{q-2}\nabla u\nabla u+\mathcal{V}(x)|u|^{q-2}uu\Big)dx\nonumber\\
  &\geq &  k_{1}\left(\int_{\mathbb{R}^{N}}\frac{1}{p}\Big(|\nabla u|^{p}+\mathcal{V}(x)|u|^{p}\Big)dx\right)^{\sigma_{1}-1}\|u\|_{\mathcal{V},p} +k_{2}\left(\int_{\mathbb{R}^{N}}\frac{1}{q}\Big(|\nabla u|^{q}+\mathcal{V}(x)|u|^{q}\Big)dx\right)^{\sigma_{2}-1}\|u\|_{\mathcal{V},q}\nonumber\\
  &\geq & \frac{k_{1}}{p}\|u\|_{\mathcal{V},p}^{p\sigma_{1}}+\frac{k_{2}}{q}\|u\|_{\mathcal{V},q}^{q\sigma_{2}}\nonumber\\
  &\geq & \frac{\min\left(k_{1},k_{2}\right)}{q}(\|u\|_{\mathcal{V},p}^{p\min(\sigma_{1},\sigma_{2})}+\|u\|_{\mathcal{V},q}^{p\min(\sigma_{1},\sigma_{2})})\nonumber\\
  &\geq & \frac{\min\left(k_{1},k_{2}\right)}{q2^{p\min(\sigma_{1},\sigma_{2})-1}}\vertiii{u}_{\mathcal{V},p,q}^{p\min(\sigma_{1},\sigma_{2})-1}.
\end{eqnarray}
Hence, $\Psi_{1}'$ is coercive.

Using Lemma \ref{alem3.1}, we can establish that $\Psi_{1}$ belongs to $C^{1}(W_{\mathcal{V}}^{p,q}(\mathbb{R}^{N}), \mathbb{R})$. Consequently, $\Psi_{1}'$ exhibits hemicontinuity. Therefore, in accordance with the Minty-Browder Theorem, we deduce the existence of ${\Psi_{1}'}^{-1}: \Big(W_{\mathcal{V}}^{p,q}(\mathbb{R}^{N})\Big)^{*}\longrightarrow W_{\mathcal{V}}^{p,q}(\mathbb{R}^{N})$, and it is subject to boundedness.

To establish the continuity of ${\Psi_{1}'}^{-1}$, we aim to demonstrate its sequential continuity. Consider a sequence $(u_n)_{n}$ in $W_{\mathcal{V}}^{p,q}(\mathbb{R}^{N})^{*}$ that strongly converges to $u \in \Big(W_{\mathcal{V}}^{p,q}(\mathbb{R}^{N})\Big)^{*}$. Consider $v_{n} = {\Psi_{1}'}^{-1}(u_{n})$, and $v = {\Psi_{1}'}^{-1}(u)$. We can establish that the sequence $(v_n)_{n}$ is bounded within $W_{\mathcal{V}}^{p,q}(\mathbb{R}^{N})$, so, we may make the assumption that it converges weakly to a certain $v_{0} \in W_{\mathcal{V}}^{p,q}(\mathbb{R}^{N})$. As a consequence of the strong convergence of $u_{n}$ to $u$, we can establish:
\begin{equation*}
    \lim_{n\rightarrow\infty}\langle \Psi_{1}'(v_{n}), v_{n}-v_{0}\rangle=\langle u_{n}, v_{n}-v_{0}\rangle=0,
\end{equation*}
that is,
\begin{equation}\label{ATA11}
    \lim_{n\rightarrow\infty} K_{1}\left(\frac{\|v_{n}\|_{\mathcal{V},p}^{p}}{p}\right)\langle \varphi_{p,\mathcal{V}}'(v_{n}), v_{n}-v_{0}\rangle=0,
\end{equation}
and
\begin{equation}\label{ATA11bis}
    \lim_{n\rightarrow\infty} K_{2}\left(\frac{\|v_{n}\|_{\mathcal{V},q}^{q}}{q}\right)\langle \varphi_{q,\mathcal{V}}'(v_{n}), v_{n}-v_{0}\rangle=0.
\end{equation}
Because ${v_{n}}$ is bounded in $W_{\mathcal{V}}^{p,q}(\mathbb{R}^{N})$, we can deduce that $\varphi_{p,\mathcal{V}}(v_{n})$ and $\varphi_{q,\mathcal{V}}(v_{n})$ are also bounded, therefore
\begin{equation*}
    \|u\|_{\mathcal{V},p}\rightarrow t_{0}\geq0,\quad \|u\|_{\mathcal{V},q}\rightarrow t_{1}\geq0, \textrm{ as } n\rightarrow\infty.
\end{equation*}
If $t_{0}= t_{1}= 0$, hence,  we establish that $(v_n)_{n}$ exhibits strong convergence to $v_{0}$ within the function space $W_{\mathcal{V}}^{p,q}(\mathbb{R}^{N})$. Thanks to the continuity and injectivity of ${\Psi_{1}'}^{-1}$. As a result, we obtain the desired conclusion.

If $t_{0}>0$ (respectively $t_{1}>0$), it follows from the continuity of the function $K_{1}$ (respectively $K_{2}$) that
\begin{equation*}
    K_{1}\left(\frac{\|v_{n}\|_{\mathcal{V},p}^{p}}{p}\right)\rightarrow K_{1}\left(\frac{t_{0}}{p}\right), \textrm{ respectively } K_{2}\left(\frac{\|v_{n}\|_{\mathcal{V},q}^{q}}{q}\right)\rightarrow K_{2}\left(\frac{t_{1}}{q}\right)\textrm{ as } n\rightarrow\infty.
\end{equation*}
Thus, by \ref{K1} (respectively \ref{K2}), when $n$ is taken to be sufficiently large, the result becomes apparent:
\begin{equation}\label{ATA12}
     K_{1}\left(\frac{\|v_{n}\|_{\mathcal{V},p}^{p}}{p}\right)\geq \tau_{0}>0, \textrm{ respectively } K_{2}\left(\frac{\|v_{n}\|_{\mathcal{V},q}^{q}}{q}\right)\geq \tau_{1}>0.
\end{equation}
As a result of \eqref{ATA11} and \eqref{ATA12}, the following is obtained:
\begin{equation}\label{ATA13}
    \lim_{n\rightarrow\infty}\langle \Psi_{1}'(v_{n}), v_{n}-v_{0}\rangle=0.
\end{equation}
Utilizing \eqref{ATA13} and considering the weak convergence of $v_{n}$ to $v_{0}$ in $W_{\mathcal{V}}^{p,q}(\mathbb{R}^{N})$, we can employ Lemma \ref{lem7} to conclude that $v_{n}$ undergoes a strong convergence to $v_{0}$ within the same function space $W_{\mathcal{V}}^{p,q}(\mathbb{R}^{N})$.
\end{prv}
\begin{center}
\bf Step 2 : Coercivity of $\Psi_{1}$
\end{center}
\begin{lem}\label{lemceorcive}
Under the hypothesis \ref{V1}, \ref{K1} and \ref{K2} the functional $\Psi_{1}$ is coercive.
\end{lem}
\begin{prv}
In cases where $\vertiii{u}_{\mathcal{V},p,q}$ is exceptionally large, it allows for the choice of $\|u\|_{\mathcal{V},p}\geq1$ and $\|u\|_{\mathcal{V},q}\geq1$. Then, by \ref{K1} and \ref{K2} we have
      \begin{eqnarray}\label{eeeeq1}
        \Psi_{1}(u) &=& \widehat{K}_{1}\left(\varphi_{p,\mathcal{V}}(u)\right)+\widehat{K}_{2}\left(\varphi_{q,\mathcal{V}}(u)\right) \nonumber\\
        &\geq& \frac{k_{1}}{p\sigma_{1}}\|u\|_{\mathcal{V},p}^{p\sigma_{1}}+ \frac{k_{2}}{q\sigma_{2}}\|u\|_{\mathcal{V},q}^{q\sigma_{2}}\nonumber\\
        &\geq& \frac{\min(k_{1},k_{2})}{q\max(\sigma_{1},\sigma_{2})}\Big(\|u\|_{\mathcal{V},p}^{p\min(\sigma_{1},\sigma_{2})}
        +\|u\|_{\mathcal{V},q}^{p\min(\sigma_{1},\sigma_{2})}\Big)\nonumber\\
        &\geq& \frac{\min(k_{1},k_{2})}{q\max(\sigma_{1},\sigma_{2})2^{p\min(\sigma_{1},\sigma_{2})-1}}\vertiii{u}_{\mathcal{V},p,q}^{p\min(\sigma_{1},\sigma_{2})}.
      \end{eqnarray}
Which shows that $\Psi_{1}$ is coercive.
\end{prv}
\begin{center}
\bf Step 3 : Application of Theorem \ref{ricceri}
\end{center}
Taking $X=W_{\mathcal{V}}^{p,q}(\mathbb{R}^{N})$, $\Psi_{1}$ and $\Psi_{2}$ are in accordance with the definitions set earlier, thanks to Lemma \ref{alem3.1} $\Psi_{2}$ is $C^{1}$-functional with compact derivative.  Additionally, Lemma \ref{alem3.2} highlights the sequential weak lower continuity of $\Psi_{1}$ and by Lemma \ref{ii} $\Psi_{1}$ is $C^{1}$-functional belonging
to $\mathcal{W}_{W_{\mathcal{V}}^{p,q}(\mathbb{R}^{N})}$. Moreover, as implied by Lemma \ref{invertible}, it's essential to note that the operator $\Psi_{1}'$ possesses a continuous inverse on $\Big(W_{\mathcal{V}}^{p,q}(\mathbb{R}^{N})\Big)^{*}$, and by Lemma \ref{lemceorcive} $\Psi_{1}$ is coercive.

It is evident that $u_{0}=0$ is global minimum of $\Psi_{1}$ and that $\Psi_{1}(u_{0})=\Psi_{2}(u_{0})=0$.
Furthermore, it is essential to note that $\Psi_{1}$ remains bounded when restricted to any bounded subset of $W_{\mathcal{V}}^{p,q}(\mathbb{R}^{N})$. In fact, if $\vertiii{u}_{\mathcal{V},p,q} \leq c$, we have:
\begin{eqnarray}\label{013}
  \Psi_{1}(u) &=& \widehat{K}_{1}\left(\varphi_{p,\mathcal{V}}(u)\right)+\widehat{K}_{2}\left(\varphi_{q,\mathcal{V}}(u)\right) \nonumber\\
  &\leq& \widehat{K}_{1}\left(\frac{\|u\|_{\mathcal{V},p}^{p}}{p}\right)+\widehat{K}_{2}\left(\frac{\|u\|_{\mathcal{V},q}^{q}}{q}\right) \nonumber\\
  &\leq& \widehat{K}_{1}\left(\frac{c^{p}}{p}\right)+\widehat{K}_{2}\left(\frac{c^{q}}{q}\right).
\end{eqnarray}
Now, in accordance with the assumptions mentioned in \ref{F2}, for any $\varepsilon > 0$, there exists a positive value $\zeta_{1}$ satisfying:
\begin{equation}\label{eq14}
    |F_{1}(x,t)|\leq \varepsilon|t|^{q\max(\sigma_{1} ,\sigma_{2} )},
\end{equation}
for every $x \in \mathbb{R}^{N}$ and $t \in [-\zeta_{1} , \zeta_{1}]$. Given that $f_{1}\in \mathcal{A}$, for $r \in[q\max(\sigma_{1} ,\sigma_{2} ), q^{*}]$, we can choose $c_{1}>0$ so that
\begin{equation}\label{eq15}
    |F_{1}(x,t)|\leq c_{1}|t|^{r},
\end{equation}
 for all  $(x, t) \in\mathbb{R}^{N}\times(\mathbb{R}\setminus[-\zeta_{1} , \zeta_{1}])$. Combining \eqref{eq14} and \eqref{eq15}, we get
\begin{equation}\label{eqeq15}
    |F_{1}(x,t)|\leq \varepsilon|t|^{q\max(\sigma_{1} ,\sigma_{2} )}+c_{1}|t|^{r},
\end{equation}
for each $(x, t)\in \mathbb{R}^{N} \times \mathbb{R}$. So by Lemma \ref{embeddingr}, the embedding $W_{\mathcal{V}}^{p,q}(\mathbb{R}^{N})$ in $L^{q\max(\sigma_{1} ,\sigma_{2} )}(\mathbb{R}^{N} )$ and $L^{r}(\mathbb{R}^{N} )$
is compact. Then for some positive constants $c_{2},$ and $c_{3}$, one has for all $u\in W_{\mathcal{V}}^{p,q}(\mathbb{R}^{N})$
\begin{eqnarray*}
  \Psi_{2}(u) &\leq& \varepsilon \|u\|^{q\max(\sigma_{1} ,\sigma_{2} )}_{L^{\max(p\sigma_{1} ,q\sigma_{2} )}(\mathbb{R}^{N})}+c_{1} \|u\|^{r}_{L^{r}(\mathbb{R}^{N})}\\
  &\leq& \varepsilon c_{2}\vertiii{u}_{\mathcal{V},p,q}^{q\max(\sigma_{1} ,\sigma_{2} )}+c_{1}c_{3}\vertiii{u}_{\mathcal{V},p,q}^{r}.
\end{eqnarray*}
Or by \ref{K1} and \ref{K2}, if $\vertiii{u}_{\mathcal{V},p,q}\leq1$, it permits the selection of values for $\|u\|_{\mathcal{V},p}<1$ and $\|u\|_{\mathcal{V},p}<1$. Then, as implied by \ref{K1} and \ref{K2}, we have:
 \begin{eqnarray}\label{eeeeq2}
 \Psi_{1}(u) &\geq& \frac{\min(k_{1},k_{2})}{q\max(\sigma_{1},\sigma_{2})}\Big(\|u\|_{\mathcal{V},p}^{p\max(\sigma_{1},\sigma_{2})}
        +\|u\|_{\mathcal{V},q}^{q\max(\sigma_{1},\sigma_{2})}\Big)\nonumber\\
 &\geq& \frac{\min(k_{1},k_{2})}{q\max(\sigma_{1},\sigma_{2})2^{q\max(\sigma_{1},\sigma_{2})-1}}\vertiii{u}_{\mathcal{V},p,q}^{q\max(\sigma_{1},\sigma_{2})}.
 \end{eqnarray}
Hance
\begin{equation*}
    \vertiii{u}_{\mathcal{V},p,q}^{q\max(\sigma_{1} ,\sigma_{2} )}\leq \frac{q\max(\sigma_{1},\sigma_{2})2^{q\max(\sigma_{1},\sigma_{2})-1}}{\min(k_{1},k_{2})}\Psi_{1}(u),
\end{equation*}
then
\begin{equation*}
    \Psi_{2}(u)\leq \varepsilon c_{2}\frac{q\max(\sigma_{1},\sigma_{2})2^{q\max(\sigma_{1},\sigma_{2})-1}}{\min(k_{1},k_{2})}\Psi_{1}(u)
    +c_{1}c_{3}\left(\frac{q\max(\sigma_{1},\sigma_{2})2^{q\max(\sigma_{1},\sigma_{2})-1}}{\min(k_{1},k_{2})}\Psi_{1}(u)\right)^{\frac{r}{q\max(\sigma_{1} ,\sigma_{2})}}.
\end{equation*}
Consequently, since $r>q\max(\sigma_{1},\sigma_{2})$ and given that $\Psi_{1}(u)\rightarrow 0$ as $u \rightarrow0$, we can deduce that:
\begin{equation}\label{eq16}
    \limsup_{u\rightarrow0}\frac{\Psi_{2}(u)}{\Psi_{1}(u)}\leq  \varepsilon c_{2}\frac{q\max(\sigma_{1},\sigma_{2})2^{q\max(\sigma_{1},\sigma_{2})-1}}{\min(k_{1},k_{2})}.
\end{equation}
On the flip side, as implied by \ref{F3}, it follows that, for all $\varepsilon > 0$, there exists a $\zeta_{2} > 0$ satisfying:
\begin{equation}\label{eq17}
    |F_{1}(x,t)|\leq \varepsilon |t|^{p\min(\sigma_{1},\sigma_{2})},
\end{equation}
for any $x\in\mathbb{R}^{N}$ and $|t| > \zeta_{2}$, and given that $f \in \mathcal{A}$, there exists a constant $c_{4} > 0$ for $r \in [p, p\min(\sigma_{1},\sigma_{2})]$, satisfying:
\begin{equation}\label{eq18}
     |F_{1}(x,t)|\leq c_{4}|t|^{r},
\end{equation}
for each $x\in\mathbb{R}^{N}$ and $|t| \leq \zeta_{2}$. If $\vertiii{u}_{\mathcal{V},p,q}>1$, we can apply Lemma \ref{lemceorcive}, resulting in:
\begin{eqnarray}\label{eq18eq19}
  \frac{\Psi_{2}(u)}{\Psi_{1}(u)} &=& \frac{\displaystyle\int_{\Omega}F_{1}(x,u)dx}{\widehat{K}_{1}\left(\varphi_{p,\mathcal{V}}(u)\right)+\widehat{K}_{2}\left(\varphi_{q,\mathcal{V}}(u)\right)} \nonumber\\
  &\leq& \frac{\displaystyle\int_{\{x\in \mathbb{R}^{N}:|u(x)|\leq\zeta_{2} \}}F_{1}(x,u)dx}{\frac{q\max(\sigma_{1},\sigma_{2})2^{p\min(\sigma_{1},\sigma_{2})-1}}{\min(k_{1},k_{2})}\vertiii{u}_{\mathcal{V},p,q}^{p\min(\sigma_{1},\sigma_{2})}} +\frac{\displaystyle\int_{\{ x\in \mathbb{R}^{N}:|u(x)|>\zeta_{2}\}}F_{1}(x,u)dx}{\frac{q\max(\sigma_{1},\sigma_{2})2^{p\min(\sigma_{1},\sigma_{2})-1}}
  {\min(k_{1},k_{2})}\vertiii{u}_{\mathcal{V},p,q}^{p\min(\sigma_{1},\sigma_{2})}}\nonumber\\
  &\leq&\frac{c_{4}\|u\|_{L^{r}(\mathbb{R}^{N})}^{r}}{\frac{q\max(\sigma_{1},\sigma_{2})2^{p\min(\sigma_{1},\sigma_{2})-1}}{\min(k_{1},k_{2})}\vertiii{u}_{\mathcal{V},p,q}^{p\min(\sigma_{1},\sigma_{2})}} +\frac{\varepsilon \vertiii{u}_{\mathcal{V},p,q}^{p\min(\sigma_{1},\sigma_{2})}}{\frac{q\max(\sigma_{1},\sigma_{2})2^{p\min(\sigma_{1},\sigma_{2})-1}}
  {\min(k_{1},k_{2})}\vertiii{u}_{\mathcal{V},p,q}^{p\min(\sigma_{1},\sigma_{2})}}\nonumber\\
  &\leq&\frac{c_{4}\|u\|_{L^{r}(\mathbb{R}^{N})}^{r}}{\frac{q\max(\sigma_{1},\sigma_{2})2^{p\min(\sigma_{1},\sigma_{2})-1}}{\min(k_{1},k_{2})}\vertiii{u}_{\mathcal{V},p,q}^{p\min(\sigma_{1},\sigma_{2})}} +\varepsilon c_{5}.
\end{eqnarray}
Since $p\min(\sigma_{1},\sigma_{2}) > r$, then,
\begin{equation}\label{eq19}
    \limsup_{\vertiii{u}_{\mathcal{V},p,q}\rightarrow\infty}\frac{\Psi_{2}(u)}{\Psi_{1}(u)}\leq \varepsilon c_{5}.
\end{equation}
As $\varepsilon > 0$ is arbitrary, the expressions \eqref{eq16} and \eqref{eq19} signify that:
\begin{equation*}
    \max \left\{\limsup_{\vertiii{u}_{\mathcal{V},p,q}\rightarrow\infty}\frac{\Psi_{2}(u)}{\Psi_{1}(u)}, \limsup_{\vertiii{u}_{\mathcal{V},p,q}\rightarrow0}\frac{\Psi_{2}(u)}{\Psi_{1}(u)} \right\}\leq0.
\end{equation*}
Thus, all the conditions of Theorem \ref{ricceri} are verified. Hence, for every compact interval $\mathbf{\Lambda} \subset (\Theta^{*}, +\infty)$, there exists a positive $\eta$ that meets the criteria outlined in the conclusion of Theorem \ref{ricceri}. Let $\alpha\in \mathbf{\Lambda}$ and $f_{2}\in \mathcal{A}$. We can set:
\begin{equation*}
    \Psi_{3}(u)=\int_{\mathbb{R}^{N}}F_{2}(x, u)dx,\quad F_{2}(x, t)=\int_{0}^{t}f_{2}(x, s)ds,
\end{equation*}
for every $u \in W_{\mathcal{V}}^{p,q}(\mathbb{R}^{N})$, we have $\Psi_{3}$ as a $C^{1}$ functional on $W_{\mathcal{V}}^{p,q}(\mathbb{R}^{N})$ with a derivative that is compact. As a result, there exists a positive $\beta^{*}$ such that, for any $\beta$ within the interval $[0, \beta^{*}]$, the equation:
\begin{equation*}
	\Psi_{1}'(u)=\alpha \Psi_{2}'(u)+\beta\Psi_{3}'(u),
	\end{equation*}
consists of at least three solutions, with the norm of each being smaller than $\eta$. Importantly, the solutions in $W_{\mathcal{V}}^{p,q}(\mathbb{R}^{N})$ for the equation above are precisely the weak solutions of problem \eqref{p}. Thus, we have completed the proof of Theorem \ref{Theorem3}.
\end{th1}

Here, we offer an illustrative example to expound upon the conclusions laid out in Theorem \ref{Theorem3}.
\begin{cor}
Let $p\in [2, N[$, $q\in [p+1, N[$ and $\xi>\max\Big(2,\max (\sigma_{1}q,\sigma_{2}q)\Big)$ such that $p<q<\xi$, where $1<\sigma_{1}<\frac{N}{N-q}$ and $1<\sigma_{2}<\frac{N}{N-q}$.

 We consider
\begin{equation}\label{20.20}
    K_{1}(t)=(1+t)^{\sigma_{1}-1}, \quad K_{2}(t)=(1+t)^{\sigma_{2}-1} \mbox{ for all } t\geq0.
\end{equation}
and
\begin{equation}\label{20.20.20}
    f_{1}(t)=\xi \sin(t)\cos(t)|\sin(t)|^{\xi-2} \mbox{ for every } t\in \mathbb{R}.
\end{equation}
Therefore, by referencing both \eqref{20.20} and \eqref{20.20.20}, we obtain the following:
\begin{equation*}
    \widehat{K}_{1}(t)=\frac{(1+t)^{\sigma_{1}}}{\sigma_{1}}, \quad \widehat{K}_{2}(t)=\frac{(1+t)^{\sigma_{2}}}{\sigma_{2}}
    \mbox{ and } F_{1}(x,t)=F_{1}(t)=|\sin(t)|^{\xi}.
\end{equation*}
Then, the following nonlinear perturbed $(p, q)$-Schr\"{o}dinger-Kirchhoff problem
\begin{eqnarray}\label{exexex}
 &&\left(1+ \left(\int_{\mathbb{R}^N}\frac{1}{p}|\nabla u|^p dx+ \frac{1}{p}\int_{\mathbb{R}^N} \mathcal{V}(x)|u|^p  dx\right)^{\sigma_{1}-1}\right)\Big(-\Delta_p u+\mathcal{V}(x)|u|^{p-2}u\Big)\nonumber\\
 &&+\left(1+ \left(\int_{\mathbb{R}^N}\frac{1}{q}|\nabla u|^q dx+ \frac{1}{q}\int_{\mathbb{R}^N} \mathcal{V}(x)|u|^q  dx\right)^{\sigma_{2}-1}\right)\Big(-\Delta_q u+\mathcal{V}(x)|u|^{q-2}u\Big)\nonumber\\
&&\quad  = \alpha \xi \sin(u)\cos(u)|\sin(u)|^{\xi-2}+\beta f_{2}(x,u) \quad \textrm{in }\mathbb{R}^N.
\end{eqnarray}
Possesses a minimum of three weak solutions, all with norms below $\eta$.
\end{cor}
\begin{prv}
For any $t\in \mathbb{R}$, we can assert that $f_{1}$ is a member of $\mathcal{A}$. This assertion is supported by the following reasoning:
\begin{equation*}
   \sup_{t\in\mathbb{R}}\frac{|f_{1}(t)|}{|t|^{s-1}}<\infty,
\end{equation*}
is true for every $p<q< s < \xi$. Conversely, the following holds:
\begin{equation*}
\limsup_{t\rightarrow0}\frac{|\sin(t)|^{\xi}}{|t|^{q\max(\sigma_{1},\sigma_{2})}}=0
\quad \textrm{and}\quad  \limsup_{t\rightarrow\infty}\frac{|\sin(t)|^{\xi}}{|t|^{p\min(\sigma_{1},\sigma_{2})}}=0.
\end{equation*}
Choose a positively measured compact set $\mathcal{U}$ in $\mathbb{R}^{N}$ and select $v \in W_{\mathcal{U}}^{p,q}(\mathbb{R}^{N})$ so that it satisfies $v(x) = \frac{\pi}{2}$ within $\mathcal{U}$ and $0 \leq v(x) \leq \frac{\pi}{2}$ outside of $\mathbb{R}^{N}\backslash \mathcal{U}$. This results in:
$$ \int_{\mathbb{R}^{N}}|\sin(v(x))|^{\xi}dx=\text{meas}(\mathcal{U})+\int_{\mathbb{R}^{N}\backslash \mathcal{U}}|\sin(v(x))|^{\xi}dx>0,$$
signifying that \ref{F1}, \ref{F2} and \ref{F3} are verified. Also, for $k_{1}=1$ and $k_{2}=1$ the condition
\ref{K1} and \ref{K2} is verified, we set
\begin{equation*}
    \Theta^{*}=\inf\left\{\frac{\frac{\left(1+\int_{\mathbb{R}^{N}}\frac{1}{p}(|\nabla v|^{p}+\mathcal{V}(x)|v|^{p})dx\right)^{\sigma_{1}}}{\sigma_{1}}
    +\frac{\left(1+\int_{\mathbb{R}^{N}}\frac{1}{q}(|\nabla v|^{q}+\mathcal{V}(x)|v|^{q})dx\right)
    ^{\sigma_{2}}}{\sigma_{2}}}{\displaystyle\int_{\mathbb{R}^{N}}|\sin(v(x))|^{\xi}dx}
    :v\in W_{\mathcal{V}}^{p,q}(\mathbb{R}^{N}), \int_{\mathbb{R}^{N}}|\sin(v(x))|^{\xi}dx>0\right\}.
\end{equation*}
Then, it follows from Theorem \ref{Theorem3}, we deduce that, within any compact interval $\mathbf{\Lambda}\subset(\Theta^{*}, \infty)$, there exist specific values $\eta>0$ and $\beta^{*}>0$, such that for each $\alpha\in\mathbf{\Lambda}$, and every $\beta\in[0, \beta^{*}]$, and for any $f_{2}\in \mathcal{A}$, there exist a minimum of three weak solutions to problem \eqref{exexex}, each having norms smaller than $\eta$.
\end{prv}
\bibliographystyle{plain}

\end{document}